\documentclass{article}
 \usepackage{amssymb}
\usepackage{amsmath}
\newcommand{\meanint}{{\int{\mkern-19mu}-}}
\def\ring{\mathaccent"0017}

\def\ring{\mathaccent"0017}

\title{Estimates for differential operators of  vector analysis involving  $L^1$-norm}
\author{Vladimir Maz'ya\footnote{The author was partially supported by the USA National Science Foundation grant DMS 0500029 and by the UK Engineering and Physical Sciences Research Council grant EP/F005563/1.}
\\*[10pt]
\emph{\small Department of Mathematical Sciences,   University of Liverpool, Liverpool L69 7ZL}
\\
\emph{\small and}\\
\emph{\small Department of Mathematics, Link\"oping University, Link\"oping, SE-581 83}
\\
\emph{\small e-mail: vlmaz@mai.liu.se}}
\date{}

\begin{document}
\maketitle

\begingroup

\narrower\noindent
{\it Abstract.}  New Hardy and Sobolev type inequalities involving $L^1$-norms of scalar and vector-valued functions in $\Bbb{R}^n$ are obtained. The work is related  to some problems stated in the recent paper by Bourgain and Brezis \cite{[BB2]}.
\endgroup

\bigskip\noindent
{ {\small\it  Mathematics Subject Classification (2000)}: }{42B20, 42B25, 46E35}

\medskip\noindent
{ {\small\it Key words}: }{Newtonian potential, Hardy-type inequality, divergence free fields, ${\rm div}-{\rm curl}$ inequalities

\section{Introduction}

Starting with the pioneering paper by Bourgain and Brezis \cite{[BB1]}, much interest arose in  various 
 $L^1$-estimates for vector fields (see \cite{[BB2]}, \cite{[BB3]}, \cite{[BV1]}, \cite{[BV2]}, \cite{[LS]}, \cite{[VS1]}--\cite{[VS4]}, \cite{[Ma2]}, \cite{[MS]} {\it et al}).  The present article belongs to the same area and it was inspired by a question  Ha\"im Brezis asked me at a recent conference in Rome. The question concerns the validity of the Hardy-type inequality
 \begin{equation}\label{00}
 \int_{\Bbb{R}^n} | D{\bf u} (x) |\frac{dx}{|x|} \leq const. \int_{\Bbb{R}^n} |\Delta {\bf u} (x)|\, dx
 \end{equation}
 in the case of divergence free $\Delta {\bf u}$ and, in a modified form, is included in {\it Open Problem} $1$ formulated in \cite{[BB2]} on p. 295.
 
 \smallskip
 
 In this paper  a positive answer to Brezis' question is given (Theorem 2)  and some  related results are obtained. For instance,  by  Theorem 1, the inequality
 \begin{equation}\label{2w}
\Bigl |\int_{\Bbb{R}^n} \sum_{1\leq j\leq n} a_j \Bigl |\frac{\partial u}{\partial x_j}\Bigr | \, \frac{dx}{|x|} \, \Bigr | \leq const. 
\int_{\Bbb{R}^n} |\Delta u|\, dx,
\end{equation}
where $a_j$ are real constants, holds for all real valued scalar functions $u\in C_0^\infty$ if and only if $$
\sum_{1\leq j\leq n} a_j  =0 .$$
 At the end of the paper certain inequalities for vector valued functions involving Hilbert-Sobolev spaces ${\cal H}^{-s}(\Bbb{R}^n)$ of negative order are collected. For example,  by Theorem 3 $(iii)$, the estimate
 \begin{equation}\label{3w}
\bigl | \|{\bf g}\|^2_{{\cal H}^{-n/2}} - n\|{\rm div} \, {\bf g} \|^2_{{\cal H}^{-1-n/2}} \bigr | \leq \frac{(2\sqrt{\pi})^{-n}}{\Gamma(n/2)} \|{\bf g}\|^2_{L^1}
\end{equation}
 holds for all ${\bf g}\in L^1$   with ${\rm div}\,  {\bf g}\in {\cal H}^{-1-n/2}$.   An assertion dual to (\ref{3w}) replies in affirmative to {\it Open Problem} $2$ on p. 297 in \cite{[BB2]} for the particular case $l=1$, $p=2$, $s=n/2$.
 
 \smallskip
 
 We make no difference in notations between spaces of scalar and vector-valued  functions. If the domain of integration is not indicated, the integral is taken over $\Bbb{R}^n$. We never mention $\Bbb{R}^n$ in notations of function spaces.

\section{Inequality for scalar functions}

\medskip	

{\bf Theorem 1.}
{\it Let $f$ and $\Phi$ denote  scalar real-valued functions defined on $\Bbb{R}^n$. Assume that $f\in L^1$ and 
\begin{equation}\label{1}
\int  f(x) dx =0.
\end{equation}
Furthermore, let $\Phi$ be Lipschitz on the unit sphere $S^{n-1}$  and positively homogeneous of degree $q\in [1, \frac{n}{n-1})$.  By $u$ we mean the Newtonian (logarithmic for $n=2$) potential of $f$:
$${ u}(x) = \int  \Gamma (x-y)\, { f}(y)\, dy,
$$
where $\Gamma(x)$ is the fundamental solution of $-\Delta$.

\smallskip

A necessary and sufficient condition for the inequality
\begin{equation}\label{2}
\sup\limits_{R>0} \Bigl |\,  \int_{|x|<R} \Phi \bigl(\nabla u (x)\bigr) |x|^{n(q-1) -q} dx\, \Bigr| \leq C \Bigl (\int  |f(x)| dx\Bigr)^q
\end{equation}
to hold for all $f$ is
\begin{equation}\label{3}
\int_{S^{n-1}} \Phi (x) \, d\omega_x =0.
\end{equation}
The constant $C$ in $(\ref{2})$ depends only on $\Phi$, $q$,  and} $n$.

\smallskip

Here and elsewhere $d\omega_x$ is the area element of the unit sphere $S^{n-1}$ at the point $x/|x|$.

\smallskip

{\bf Conjecture.} It seems plausible that the inequality (\ref{2}) holds also for the critical value $q= n/(n-1)$.  The following simple assertion obtained in \cite{[MS]} speaks in favour. The inequality
$$\Bigl |\int_{\Bbb{R}^2} \sum _{i,j =1}^2 a_{ij}\frac{\partial u}{\partial x_i} \frac{\partial {u}}{\partial x_j}\, dx \Bigr | \leq C \Bigl(\int_{\Bbb{R}^2}|\Delta u |\, dx \Bigr)^2$$
with $a_{ij} = const$ holds for all $u\in C_0^\infty$ if and only if $a_{11} + a_{22} =0$.

\smallskip

{\bf Proof of Theorem 1.} The necessity of (\ref{3}) can be  derived  by putting a sequence of radial mollifications of the Dirac function in place of $f$ in (\ref{2}).

Let us  prove the sufficiency of (\ref{3}).  We write $\nabla u (x)$ in the form
$$\nabla u (x) = \sum_{j=1}^4 A_j(x),$$
where
\begin{eqnarray*}
A_1(x)& =&  \frac{1}{|S^{n-1}|} \int_{|y|<|x|/2}\Bigl( \frac{y-x}{|y-x|^n} + \frac{x}{|x|^n}\Bigr) f(y)\, dy,\\
\\
A_2(x) &= &\frac{1}{|S^{n-1}|} \int_{|x|/2<|y|<2|x|} \frac{y-x}{|y-x|^n}  f(y)\, dy,\\
\\
A_3(x) &=& \frac{1}{|S^{n-1}|} \int_{|y|>2|x|} \frac{y-x}{|y-x|^n}  f(y)\, dy,
\end{eqnarray*}
and 
$$A_4(x) = \frac{-1}{|S^{n-1}|}\frac{x}{|x|^n}  \int_{|y|<|x|/2} f(y)\, dy.$$
By (\ref{3}), for all $R>0$,
\begin{equation}\label{4}
\int_{|x|<R} \Phi\bigl( A_4(x)\bigr) |x|^{n(q-1) -q} dx =0.
\end{equation}
We check directly that
\begin{eqnarray}
|A_1(x)| & \leq & \frac{c}{|x|^n}\int_{|y|<|x|/2} |f(y)|\, |y|\, dy,\label{n1}\\
\nonumber\\
|A_2(x)| & \leq & c \int_{|x|/2<|y|<2|x|}  \frac{|f(y)|}{|y-x|^{n-1}} dy,\label{n2}\\
\nonumber\\
|A_3(x)| & \leq & c \int_{|y|>2|x|}  |f(y)| \frac{dy}{|y|^{n-1}}\label{n3}.
\end{eqnarray}
(Here and elsewhere, by $c$ we denote constants depending only on $n$ and $q$.) 
Hence
\begin{eqnarray}
&&\int  \sum _{j=1}^3 |A_j(x)| \frac{dx}{|x|}\nonumber\\
\nonumber\\
&& \leq c \int  |f(y)|\Bigl( |y|\int_{|x|>2|y|} \frac{dx}{|x|^{n+1}} + \int_{|y|/2<|x|<2|y|}\frac{dx}{|x| |x-y|^{n-1}} \nonumber\\
\nonumber\\
&&+ \frac{1}{|y|^{n-1}}\int_{|x|<|y|/2}\frac{dx}{|x|}\Bigr) dy
  \leq c \int |f(y)| \, dy\label{n4}
\end{eqnarray}

 Since $\Phi$ is Lipschitz on $S^{n-1}$ and positively homogeneous of degree $q$, we have 
 $$|\Phi(a+b) - \Phi (a)| \leq C_\Phi (|a|^{q-1} |b| + |b|^q)$$
for all $a$ and $b$ in $\Bbb{R}^n$. Now,  we deduce from (\ref{4}) that the left-hand side of (\ref{2}) does not exceed
\begin{equation}\label{n5}
c\, C_\Phi\Bigl( \int \sum _{j=1} ^3 |A_j(x)|\, |A_4(x)|^{q-1} |x|^{n(q-1) -q} dx + \int \sum _{j=1} ^3 |A_j(x)|^q  |x|^{n(q-1) -q} dx\Bigr).
\end{equation}
Because of (\ref{n4}), the first integral in (\ref{n5}) is dominated by
\begin{equation}\label{n6}
c\, \|f\|_{L^1}^{q-1}\int \sum _{j=1} ^3 |A_j(x)|\,\frac{dx}{|x|}\leq c\, C_\Phi \|f\|_{L^1}^q.
\end{equation}

Let us turn to the second integral in (\ref{n5}).
 We deduce from (\ref{n1}) and Minkowski's inequality that
\begin{equation}\label{6x}
\|{ A}_1\|_{L^q(|x|^{n(q-1) -q} dx)} \leq c  \int  |y|\, | { f} (y)|\Bigl( \int_{|x|>2|y|}\frac{dx}{|x|^{n+q}}\Bigr )^{1/q} dy.
\end{equation}
Similarly, by (\ref{n2})
\begin{equation}\label{7x}
\|{ A}_2\|_{L^q(|x|^{n(q-1) -q} dx)} \leq c  \int  | { f} (y)|\Bigl( \int_{2|y|>|x|>|y|/2} \frac{|x|^{n(q-1) -q} dx}{|y-x|^{(n-1)q}}\Bigr) ^{1/q} dy
\end{equation}
and by (\ref{n3})
\begin{equation}\label{8x}
\|{ A}_3\|_{L^q(|x|^{n(q-1) -q} dx)} \leq c  \int   | { f} (y)|\Bigl( \int_{|x|<|y|/2}|x|^{n(q-1) -q} dx\Bigr) ^{1/q} \frac{dy}{|y|^{n-1}}.
\end{equation}
Every  right-hand side  in (\ref{6x}) - (\ref{8x}) is majorized by $c\, \| { f}\|_{L^1}$. Therefore
\begin{equation}\label{9x}
\sum_{k=1}^3\|{ A}_k\|_{L^q(|x|^{n(q-1) -q} dx)} \leq  c\, \| { f} \|_{L^1}.
\end{equation}

\noindent
The proof is complete.

\section{Inequalities for vector functions}

We turn to a generalization of the inequality (\ref{00}).

\medskip

{\bf Theorem 2.} {\it Let ${\bf f}$ be an $n$-dimensional vector-valued  function in $L^1$ subject to
\begin{equation}\label{1x}
{\rm div}\,  {\bf f} =0.
\end{equation}
Also, let ${\bf u}$ denote the solution of $-\Delta {\bf u} ={\bf f}$ in $\Bbb{R}^n$ represented in the form
$${\bf u}(x) = \int  \Gamma (x-y)\, {\bf f}(y)\, dy.
$$

Then there is a constant $c$ depending on $n$ and $q\in [1, \frac{n}{n-1})$ such that}
\begin{equation}\label{2x}
\Bigl(\int   |D {\bf u}(x) |^q |x|^{n(q-1) -q} dx\Bigr)^{1/q} \leq c \int  |{\bf f}(x)|\, dx,
\end{equation}
where $D{\bf u}$ is the Jacobi matrix  $(\partial u_i/\partial x_j)_{i,j=1}^n$.

\smallskip

{\bf Remark 1.} 
The case $q\in (1, n/(n-1))$ in Theorem 2 is a consequence of the marginal cases $q=1$ and $q= n/(n-1)$ because of the H\"older inequality
$$\|\varphi \|_{L^q(|x|^{n(q-1) -q} dx)} \leq  \|\varphi \|_{L^1(|x|^{-1}dx)}^{1-n(1-1/q)}\,  \|\varphi \|_{L^{n/(n-1)}}^{n(1-1/q)}.$$
However, we prefer to deal with  all values of $q$ on the interval $[1, n/(n-1) )$ simultaneously and independently of the deeper case $q= n/(n-1)$ treated in \cite{[BB2]}).

\medskip

{\bf Proof of Theorem 2.} It follows from ${\bf f}\in L^1$ that the Fourier transform $\hat{\bf f}$ is continuous. Since  $\xi \cdot  \hat{\bf f}(\xi)=0$ by (\ref{1x}), we have $|\xi|^{-1} \xi \cdot \hat{\bf f}(0)=0$  for all $\xi\in \Bbb{R}^n\backslash \{0\}$, which is equivalent to 
\begin{equation}\label{0}
\int {\bf  f}(y)\, dy =0.
\end{equation}
(The implication (\ref{1x}) $\Longrightarrow$ (\ref{0}) was noted in \cite{[BV1]}). 

\smallskip

By the integral representation ${\bf u}= (-\Delta)^{-1} {\bf f} $ we have 
$$\Bigl | \frac{\partial {\bf u}}{\partial x_k} (x)\Bigr | \leq  \frac{1}{|S^{n-1}|}\Bigl | \int  \frac{y_k-x_k}{|y-x|^n}  {\bf f}(y)\, dy\Bigr |.$$
Obviously,
\begin{equation}\label{4x}
\Bigl | \frac{\partial {\bf u}}{\partial x_k} (x)\Bigr | \leq  \frac{1}{|S^{n-1}|}\sum_{k=1}^4{\cal A}_k(x), 
\end{equation}
where
\begin{eqnarray*}
{\cal A}_1(x)& =&  \Bigl | \int_{|y|<|x|/2}\Bigl( \frac{y_k-x_k}{|y-x|^n} + \frac{x_k}{|x|^n}\Bigr) {\bf f}(y)\, dy\Bigr |,\\
{\cal A}_2(x) &= & \Bigl | \int_{|x|/2<|y|<2|x|} \frac{y_k-x_k}{|y-x|^n}  {\bf f}(y)\, dy\Bigr |,\\
{\cal A}_3(x) &=& \Bigl |  \int_{|y|>2|x|} \frac{y_k-x_k}{|y-x|^n}  {\bf f}(y)\, dy\Bigr |,
\end{eqnarray*}
and 
\begin{equation}\label{5x}
{\cal A}_4(x) = \frac{1}{|x|^{n-1}}\Bigl |  \int_{|y|<|x|/2} {\bf f}(y)\, dy\Bigr |.
\end{equation}

\noindent
Clearly, ${\cal A}_1$, ${\cal A}_2$ and ${\cal A}_3$ satisfy (\ref{6x})-(\ref{8x}) with $f$ replaced by ${\bf f}$. Therefore,
 by Minkowski's inequality (see the proof of (\ref{9x})), we have
\begin{equation}\label{9xx}
\sum_{k=1}^3\|{\cal A}_k\|_{L^q(|x|^{n(q-1) -q} dx)} \leq  c\, \| {\bf f} \|_{L^1}.
\end{equation}

Let the $n\times n$ skew-symmetric matrix ${\cal F}$ be defined by 
$${\cal F}:= {\rm curl}\, {\bf u}: = \Bigl( \frac{\partial u_i}{\partial x_j} - \frac{\partial u_j}{\partial x_i}\Bigr)_{i,j =1}^n,$$
 i.e., 
\begin{equation}\label{3x}
{\cal F}:  = {\rm curl}\, (-\Delta)^{-1} {\bf f},
\end{equation}
where $(-\Delta)^{-1}$ stands for the Newtonian  (logarithmic for $n=2$) potential. Let ${\cal F} = (F_{ij})_{i,j =1}^n$ and ${\bf F}_j = ( F_{1j}, \ldots, F_{nj})^t$ with $t$ indicating the transposition of a matrix. We need the row divergence of the matrix ${\cal F}$:
$${\rm Div}\, {\cal F} = ({\rm div}\, {\bf F}_1, \ldots , {\rm div}\, {\bf F}_n).$$
Since

$$({\rm Div}\, {\rm curl})^t = \nabla {\rm div} - \Delta \qquad {\rm and} \qquad {\rm div}\, {\bf f} =0,$$
we have
\begin{equation}\label{3xa}
{\rm Div}\, {\cal F} = {\rm Div}\, {\rm curl}\,  (-\Delta)^{-1} {\bf f}={\bf f}^t.
\end{equation}

\smallskip

We turn to ${\cal A}_4(x)$ defined in  (\ref{5x}).  By (\ref{3xa}), we obtain from Green's formula that
\begin{equation}\label{3xb}
{\cal A}_4(x) = \frac{1}{|x|^{n-1}}\Bigl |  \int_{|y|<|x|/2} {\rm Div}\, {\cal F} (y)\, dy \Bigr | \leq c\int_{|y| = |x|/2} |{\cal F}(y)| d\omega_y,
\end{equation}
where $|{\cal F}|$ is the matrix norm. The result will follow from (\ref{9xx}), (\ref{3xb}),  and the next lemma.

\medskip

{\bf Lemma.} {\it Let ${\cal F}$ be the same skew-symmetric matrix field as  in Theorem  $2$. Then
\begin{equation}\label{10x}
\Bigl( \int  |{\cal F}(x)|^q |x|^{n(q-1) -q} dx\Bigr)^{1/q} \leq c  \int  |{\rm Div}\, {\cal F}(x)| dx,
\end{equation}
where $q\in [1, n/(n-1))$ and $c$ depends only on $n$ and $q$.}

\smallskip

{\bf Proof.} Using (\ref{3x}) and  (\ref{3xa}),  we have
\begin{eqnarray}\label{11x}
{\cal F}(x) & = & \bigl ( {\rm curl} (-\Delta)^{-1} ({\rm Div}\, {\cal F})^t \bigr ) (x)\nonumber\\
\nonumber\\
&=& \Bigl(\int_{E_1} + \int_{E_2} + \int_{E_3}\Bigr) {\rm curl}_x \bigl( (\Gamma(x-y) - \Gamma(x)\bigr) ({\rm Div}\, {\cal F}(y))^t dy,
\end{eqnarray}
where
$$E_1 = \{y: |y|\leq |x|/2\}, \,\, E_2 = \{y: |x|/2< |y|< 2|x| \}, \,\, {\rm and} \,\, E_3 = \{y: |y|\geq 2|x|\}.$$
Obviously, the norm of the part of the matrix integral (\ref{11x}) taken over $E_1$ does not exceed
\begin{equation}\label{12x}
\frac{c}{|x|^{n-1}} \int_{|y| < |x|/2} | {\rm Div}\, {\cal F} (y)| dy
\end{equation}

\noindent
and the norm of the integral over $E_2$ is dominated  by
\begin{equation}\label{13x}
c\int_{|x|/2<|y| <2|x| } | {\rm Div}\, {\cal F} (y)| \frac{dy}{|x-y|^{n-1}}.
\end{equation}

We write the part of the integral (\ref{11x}) extended over $E_3$ as
\begin{equation}\label{14x}
\int_{E_3}{\rm curl}_x\bigl(\Gamma(x-y)({\rm Div} {\cal F}(y))^t\bigr)dy + \int_{E_3}{\rm curl}_x\bigl(-\Gamma(x)({\rm Div} {\cal F}(y))^t\bigr)dy.
\end{equation}
The matrix norm of the first term in (\ref{14x}) does not exceed
\begin{equation}\label{15x}
c\int_{|y|>2|x|} |{\rm Div}{\cal F}(y) | \frac{dy}{|y|^{n-1}}.
\end{equation}

Let us denote the second integral in (\ref{14x}) by ${\cal G}(x)$ and let us put
$${\cal G} = (\mathbf{G}_1, \ldots , \mathbf{G}_n), \qquad {\rm where} \qquad {\bf G}_j = (G_{1j}, \ldots , G_{nj})^t.$$
Estimating the $L^q(|x|^{n(q-1) -q} dx)$-norms of the majorants (\ref{12x}), (\ref{13x}), and (\ref{15x}) by Minkowski's inequality, 
in the same way as we did for ${\cal A}_1$, ${\cal A}_2$, and ${\cal A}_3$ in the proof of Theorem 2, we obtain
\begin{equation}\label{16x}
\|{\cal F} - {\cal G} \|_{L^q(|x|^{n(q-1) -q} dx)} \leq c\|{\rm div} {\cal F}\|_{L^1}.
\end{equation}

\noindent
By definitions of ${\rm curl}$ and ${\rm Div}$,
\begin{eqnarray*}
G_{ij} (x) & = & \frac{\partial \Gamma}{\partial x_j} (x) \int_{E_3} {\rm div}\,  {\bf F}_i(y)\, dy -\frac{\partial \Gamma}{\partial x_i} (x)\int_{E_3} {\rm div}\,  {\bf F}_j(y)\, dy\\
\\
&=& |S^{n-1}|^{-1} |x|^{1-n} \Bigl(\frac{x_i}{|x|}\int_{E_3} {\rm div} {\bf F}_j(y)\, dy - \frac{x_j}{|x|}\int_{E_3} {\rm div}\,  {\bf F}_j(y)\, dy\Bigr)
\end{eqnarray*}
and by Green's formula,
\begin{equation}\label{17x}
G_{ij} (x) = \frac{2^{n-1}}{|S^{n-1}|}\Bigl(\frac{x_i}{|x|}\int_{|y| = 2|x|}\Bigl(\frac{y}{|y|},  {\bf F}_j(y)\Bigr) d\omega_y - 
\frac{x_j}{|x|}\int_{|y| = 2|x|}\Bigl( \frac{y}{|y|},  {\bf F}_i(y) \Bigr)d\omega_y\Bigr),
\end{equation}

\noindent
where $(\cdot \, , \cdot)$ stands for the inner product in $\Bbb{R}^n$. 
Obviously,
\begin{eqnarray*}
\int_{|z|=|x|} G_{ij} (z)\frac{z_i}{|z|} d\omega_z &=& \frac{2^{n-1}}{|S^{n-1}|}\Bigl(|S^{n-1}| \int_{|y|=2|x|}\Bigl( \frac{y}{|y|},  {\bf F}_j(y) \Bigr)d\omega_y\\
\\
&&- \int_{S^{n-1}}\frac{z_i z_j}{|z|^2} d\omega_z\int_{|y|=2|x|}\Bigl( \frac{y}{|y|},  {\bf F}_i(y) \Bigr)d\omega_y\Bigr)
\end{eqnarray*}
and since 

$$ \int_{S^{n-1}}\frac{z_i z_j}{|z|^2} d\omega_z = \frac{\delta_i^j}{n}|S^{n-1}|  ,
$$
we obtain

\begin{equation}\label{18x}
\int_{|z|=|x|}\Bigl( \frac{z}{|z|},  {\bf G}_j(z) \Bigr) d\omega_z = 2^{n-1}\frac{n-1}{n}\int_{|y| = 2|x|}\Bigl(\frac{y}{|y|},  {\bf F}_j(y)\Bigr)  d\omega_y.
\end{equation}

\noindent
For an arbitrary $r>0$ and a vector function ${\bf v}$ we set
$${\cal P}({\bf v}; r): = \int_{|y|=r} \frac{y}{|y|} {\bf v} (y) d\omega_y.$$
Now, using the majorants (\ref{12x}), (\ref{13x}), and (\ref{15x}), we deduce  from (\ref{11x}) and the definition of ${\cal G}$ that
\begin{eqnarray*}
&&\Bigl |{\cal P}({\bf F}_j; |x|) -  {\cal P}({\bf G}_j; |x|) \Bigr |\\
\\
&&\leq  c\Bigl( \frac{1}{|x|^n}\int_{E_1}|{\rm Div}\,  {\cal F} (y) | dy + \int_{E_2}|{\rm Div} \, {\cal F} (y) | \frac{dy}{|x-y|^{n-1}} + \int_{E_3}|{\rm Div} \, {\cal F} (y) | \frac{dy}{|y|^{n-1}} \Bigr).
\end{eqnarray*}
By (\ref{18x}) the left-hand side can be written in the form
$$\Bigl |{\cal P}({\bf F}_j; |x|) - 2^{n-1}\frac{n-1}{n}{\cal P}({\bf F}_j; 2|x|) \Bigr |.$$
Using the same argument as at the end of the proof of Theorem 1, we arrive at
$$\Bigl(\int \Bigl |{\cal P}({\bf F}_j; |x|) - 2^{n-1}\frac{n-1}{n}{\cal P}({\bf F}_j; 2|x|) \Bigr |^q |x|^{n(q-1) -q} dx\Bigr)^{1/q} \leq c_0\int |{\rm div}{\cal F}(x)| dx$$
which yields
$$\Bigl|\Bigl(\int \!\bigl |{\cal P}({\bf F}_j; |x|)  \bigr |^q |x|^{n(q-1) -q} dx\Bigr)^{\!1/q} - 2^{n-1}\frac{n\!-\!1}{n}\Bigl (\int |{\cal P}({\bf F}_j; 2|x|)|^q |x|^{n(q-1) -q} dx\Bigr)^{\!1/q}\Bigr|$$
$$\leq c_0\int |{\rm div}{\cal F}(x)| dx.$$
Replacing $2x$ by $x$ in the second integral of the last inequality, we can simplify  this inequality  to the form
\begin{equation}\label{19x}
\Bigl(\int \!\bigl |{\cal P}({\bf F}_j; |x|)  \bigr |^q |x|^{n(q-1) -q} dx\Bigr)^{\!1/q} \leq n\, c_0\int |{\rm div}\, {\cal F}(x)| dx.
\end{equation}
By (\ref{17x}) and (\ref{19x}), 
$$\|{\bf G}_j\|_{L^q(|x|^{n(q-1) -q} dx)} \leq c\,  \Bigl(\int \!\bigl |{\cal P}({\bf F}_j; |x|)  \bigr |^q |x|^{n(q-1) -q} dx\Bigr)^{\!1/q} 
\leq  c\int |{\rm div}\, {\cal F}(x)| dx$$
which together with (\ref{16x}) completes the proof.

\section{Generalization of Theorem 2}

In this section we show that  Theorem 2 can be extended to the vector fields ${\bf f}$, which are not necessarily divergence free.

\smallskip

First, let us collect some notation and known facts to be used in the sequel. Let $B_R = \{x\in \Bbb{R}^n: |x| <R\}$. The mean value of the integral with respect to a finite measure will be denoted by  the integral with bar. By $\hat{\varphi}$ we denote the Fourier transform of the distribution $\varphi$ (see Sect. 7.1  in \cite{[H]}).

\smallskip

The space of distributions $\varphi$ with $\nabla\varphi \in L^1$ will be denoted by $L^1_1$. This space  is endowed with the seminorm $\|\nabla\varphi\|_{L^1}$. It is well known and can be easily proved that  the finite limit
$$\varphi_\infty: = \lim\limits_{R\to \infty}\meanint_{|x|=R} \varphi(x)d\omega_x$$
exists for every $\varphi\in L^1_1$. Furthermore, $\varphi_\infty =0$ is equivalent to the inclusion of $\varphi$ in the closure $\ring{L^1_1}$ of $C_0^\infty$ in $L^1_1$.

\smallskip

The weighted Sobolev-type inequality for  all $\varphi\in \ring{L^1_1}$
\begin{equation}\label{38p}
\|\varphi\|_{L^q(|x|^{n(q-1) -q}dx)} \leq c\,  \|\nabla \varphi\|_{L^1}
\end{equation}
with $q\in [1, \frac{n}{n-1})$ can be found,  for example, in Corollary 2.1.6  \cite{[Ma1]}.

\smallskip

We formulate and prove a  result  concerning the case $q>1$. 

\medskip

{\bf Proposition 1.} {\it Let $q\in (1, \frac{n}{n-1})$ and let ${\bf u} = (-\Delta)^{-1} {\bf f}$, where ${\bf f}$ is a  vector field in $L^1$ subject to $(\ref{0})$. Also let 
$$h: = {\rm div}\,  {\bf f} \quad {\rm and} \quad 
 \nabla(-\Delta)^{-1}h \in L^1.$$
   Then}
\begin{equation}\label{Nz}
\|D {\bf u} \|_{L^q(|x|^{n(q-1) -q}dx)} \leq c \bigl ( \|{\bf f}\|_{L^1} + \|\nabla(-\Delta)^{-1} h\|_{L^1}\bigr).
\end{equation}

\smallskip

{\bf Proof.} Note that the vector function  $-\xi |\xi|^{-2} (\hat{\bf f}(\xi), \xi)$ is the Fourier transform of $\nabla(-\Delta)^{-1} h$ and that it is equal to zero at the point $\xi =0$ since $\hat{\bf f} (0) =0$. Hence
$$\int \nabla(-\Delta)^{-1} h (y)\, dy =0.$$
We see that the vector field ${\bf f} + \nabla(-\Delta)^{-1} h$ is divergence free and  integrable. Therefore, by Theorem 2,
\begin{equation}\label{N1z}
\|D(-\Delta)^{-1} \bigl( {\bf f} + \nabla(-\Delta)^{-1} h\bigr) \|_{L^q(|x|^{n(q-1) -q}dx)} \leq c\bigl(\|{\bf f}\|_{L^1} + \|\nabla(-\Delta)^{-1} h\|_{L^1}\bigr),
\end{equation}
which  implies
\begin{eqnarray}\label{N2z}
\|D {\bf u} \|_{L^q(|x|^{n(q-1) -q}dx)}& \leq& c \Bigl ( \|{\bf f}\|_{L^1} + \|D(-\Delta)^{-1}\nabla(-\Delta)^{-1} h\|_{L^q(|x|^{n(q-1) -q}dx)} \nonumber\\
\nonumber\\
&+&\|\nabla(-\Delta)^{-1} h\|_{L^1}\Bigr).
\end{eqnarray}
Since the singular integral operator $D(-\Delta)^{-1}\nabla$ is continuous in $L^q(|x|^{n(q-1) -q}dx)$ for $q\in (1, \frac{n}{n-1})$ (see \cite{[St1]}),  we derive from  (\ref{N2z}) that
\begin{eqnarray}\label{N2zz}
\|D {\bf u} \|_{L^q(|x|^{n(q-1) -q}dx)}& \leq& c \Bigl ( \|{\bf f}\|_{L^1} + \|(-\Delta)^{-1} h\|_{L^q(|x|^{n(q-1) -q}dx)} \nonumber\\
\nonumber\\
&+&\|\nabla(-\Delta)^{-1} h\|_{L^1}\Bigr).
\end{eqnarray}
Recalling that $h= {\rm div}\, {\bf f}$, we have
\begin{eqnarray*}
&&\meanint_{B_{2R}\backslash B_R} |(-\Delta)^{-1} h(x) | dx \leq c\, R^{-n} \int |{\bf f}(y)|\int_{B_{2R}\backslash B_R} \frac{dx}{|x-y|^{n-1}} dy\\
\\
&&\leq c\Bigl(R^{1-n}\int_{B_R} |{\bf f}(y)| \, dy+ \int_{\Bbb{R}^n\backslash B_R} |{\bf f}(y)|\frac{dy}{|y|^{n-1}}\Bigr).
\end{eqnarray*}
Hence
$$\lim\limits_{R\to \infty}\meanint_{B_{2R}\backslash B_R}|(-\Delta)^{-1} h(x) | dx =0$$

\noindent
and by $(-\Delta)^{-1} h\in L^1_1$ we see that
$$\lim\limits_{R\to \infty}\meanint_{|x| =R} (-\Delta)^{-1} h(x) \,  d\omega_x =0,$$
i.e., $(-\Delta)^{-1} h\in \ring{L^1_1}$. 

\smallskip

Using (\ref{38p}),  we remove the second norm on the right-hand side of (\ref{N2zz}) by changing the value of the factor $c$. The result follows. $\square$

\medskip

We turn to  the case $q=1$  which is more technical being based on properties of the Riesz transform in the Hardy space ${\bf H}$.

\smallskip

By definition, the space ${\bf H}$ consists of all integrable functions orthogonal to $1$ and is endowed with the norm
\begin{equation}\label{au}
\|\varphi\|_{{\bf H}} = \|\varphi\|_{L^1} + \|\nabla(-\Delta)^{-1/2}\varphi\|_{L^1}. 
\end{equation}
This space can be introduced also as the completion  in the norm (\ref{au}) of the set of functions $\varphi$ such that $\hat{\varphi} \in C_0^\infty(\Bbb{R}^n \backslash \{0\})$ 
(see \cite{[St2]}, Sect. 3). 

The result concerning $q=1$ which is analogous to Proposition 1 is stated as follows.

\medskip

{\bf Proposition 2.} {\it Let ${\bf u} = (-\Delta)^{-1} {\bf f}$, where ${\bf f}$ is a vector field in $L^1$ subject to  $(\ref{0})$. Also let 
$$h: = {\rm div}\,  {\bf f} \quad {\rm and} \quad 
(-\Delta)^{-1/2}h \in {\bf H}.$$
   Then}
\begin{equation}\label{Nt}
\|D {\bf u} \|_{L^1(|x|^{-1}dx)} \leq c \bigl ( \|{\bf f}\|_{L^1} + \|(-\Delta)^{-1/2} h\|_{{\bf H}}\bigr).
\end{equation}

\smallskip

{\bf Proof.} Let us show that
\begin{equation}\label{43a}
D(-\Delta)^{-1} \nabla(-\Delta)^{-1} h \in L^1_1
\end{equation}
if $(-\Delta)^{-1/2} h \in {\bf H}$. The Fourier transform of 
$$\frac{\partial ^2}{\partial x_i\partial x_j}(-\Delta)^{-1} \frac{\partial}{\partial x_k}(-\Delta)^{-1} h$$
 equals 
 $$\xi_i\xi_j |\xi|^{-2} \xi_k |\xi|^{-2} (\hat{\bf f}(\xi), \xi)$$
 
 \noindent
  and vanishes at $\xi =0$ because  $\hat{\bf f} (0) =0$. Furthermore, by definition of ${\bf H}$ and the continuity of $\nabla(-\Delta)^{-1/2}$ in ${\bf H}$ (see \cite{[St2]}, Sect. 3.4), we obtain
\begin{eqnarray}\label{Ns}
\|\partial_{x_i} D(-\Delta)^{-1}  \nabla(-\Delta)^{-1} h\|_{L^1} &\leq & \|D(-\Delta)^{-1/2}  \nabla(-\Delta)^{-1/2} (-\Delta)^{-1/2} h\|_{{\bf H}}\nonumber\\
\nonumber\\
&\leq & c\, \|(-\Delta)^{-1/2} h\|_{{\bf H}},
\end{eqnarray}
i.e. (\ref{43a}) holds. 
Next we check that the mean value of $D(-\Delta)^{-1}  \nabla(-\Delta)^{-1}h $ on the sphere $\partial B_R$ tends to zero as $R\to \infty$.  Since $h= {\rm div}\, {\bf f}$, it follows that pointwise

$$|D(-\Delta)^{-1}  \nabla(-\Delta)^{-1}h | \leq c\, (-\Delta)^{-1/2} |{\bf f}|.$$
Therefore,
\begin{eqnarray*}
&&\meanint_{B_{2R}\backslash B_R} |D(-\Delta)^{-1} \nabla(-\Delta)^{-1} h(x) | dx \leq c\, R^{-n} \int |{\bf f}(y)|\int_{B_{2R}\backslash B_R} \frac{dx}{|x-y|^{n-1}} dy\\
\\
&&\leq c\Bigl(R^{1-n}\int_{B_R} |{\bf f}(y)| \, dy+ \int_{\Bbb{R}^n\backslash B_R} |{\bf f}(y)|\frac{dy}{|y|^{n-1}}\Bigr).
\end{eqnarray*}
Hence
$$\lim\limits_{R\to \infty}\meanint_{B_{2R}\backslash B_R}|D (-\Delta)^{-1} \nabla(-\Delta)^{-1}h(x) | dx =0,$$

\noindent
which ensures that the mean value just mentioned tends to zero. 
Now we can conclude that 
$$D (-\Delta)^{-1} \nabla(-\Delta)^{-1}h \in \ring{L^1_1}.$$
 This inclusion,  together with (\ref{38p}) for $q=1$ and (\ref{Ns}), shows that the second norm on the right-hand side of (\ref{N2z}) does not exceed

$$c\, \|D (-\Delta)^{-1} \nabla(-\Delta)^{-1}h \|_{L^1_1} \leq c_1 \, \|(-\Delta)^{-1/2}h \|_{{\bf H}}.$$
As for the third norm, it has the majorant $c_2 \, \|(-\Delta)^{-1/2}h \|_{{\bf H}}$ by definition of ${\bf H}$. The result follows. 

\section{Inequalities involving $L^2$ Sobolev norms of negative order}

In the sequel, the notation ${\cal H}^l$ will be used for the space of distributions  $h$ with  finite norm

 $$\|h\|_{{\cal H}^l} : = \Bigl(\int |{ {\hat h}} (\xi)|^2 |\xi|^{2l}\, d\xi\Bigr)^{1/2},$$
 where $l\in \Bbb{R}^1$.
  
\smallskip

By $|\cdot|$ and $(\cdot, \cdot)$ the norm and the inner product in  the complex Euclidean space  will be denoted.

\smallskip

{\bf Theorem 3.}  {\it Let ${\bf g}\in C_0^\infty$ and 
\begin{equation}\label{NN}
{\bf g}_\varepsilon(x): = {\bf g}(x) - (2\pi)^{-n/2} \varepsilon^n e^{-|\varepsilon x|^2/2} \int {\bf g}(y) \, dy.
\end{equation}
Then 

$(i)$ The following limit exists and satisfies the inequality
\begin{equation}\label{39w}
\Bigl | \lim\limits_{\varepsilon\to 0_+}\bigl(\|{\bf g}_\varepsilon\|^2_{{\cal H}^{-n/2}} - n\|{\rm div} \, {\bf g}_\varepsilon \|^2_{{\cal H}^{-1-n/2}} \bigr) \Bigr | \leq \frac{(n-1) (2\sqrt{\pi})^{-n}}{\Gamma(1+ n/2)} \|{\bf g}\|^2_{L^1}.
\end{equation}

\smallskip

$(ii)$ The inequality}
\begin{equation}\label{N1y}
 \mathop{\hbox {lim sup}}_{\varepsilon\to 0_+} 
\, \bigl | \|{\bf g}_\varepsilon\|^2_{{\cal H}^{-n/2}} - c_1\, \|{\rm div} \, {\bf g}_\varepsilon \|^2_{{\cal H}^{-1-n/2}} \bigr | \leq c_2\,  \|{\bf g}\|^2_{L^1}
\end{equation}
{\it with certain constants $c_1$ and $c_2$ implies $c_1 =n$.  The constant $c_2$ satisfies
\begin{equation}\label{N2y}
 c_2\geq \frac{(n-1) (2\sqrt{\pi})^{-n}}{\Gamma(1+ n/2)},
\end{equation}
i.e. $(\ref{39w})$ is sharp.

\smallskip

$(iii)$ If
\begin{equation}\label{Nd}
\int {\bf g}(y) \, dy =0,
\end{equation}
then}
\begin{equation}\label{46a}
\Bigl | \|{\bf g} \|^2_{{\cal H}^{-n/2}} -  n\, \|{\rm div} \, {\bf g} \|^2_{{\cal H}^{-1-n/2}} \Bigr | \leq \frac{ (2\sqrt{\pi})^{-n}}{\Gamma( n/2)} \|{\bf g}\|^2_{L^1}.
\end{equation}

\smallskip

{\bf Proof.} $(i)$  The expression in parentheses on the left-hand side of 
 (\ref{39w}) can be written as
\begin{eqnarray}\label{2y}
&&(2\pi)^{-n}\Bigl (\int |\hat{\bf g}_\varepsilon(\xi)|^2 \frac{d\xi}{|\xi|^n} -n \int |\,  (\hat{\bf g}_\varepsilon(\xi), \xi)|^2 \frac{d\xi}{|\xi|^{2+n}} \Bigr ) \nonumber\\
\nonumber\\
&&= (2\pi)^{-n}\Bigl (\sum_{1\leq j,k\leq n}\int\ \frac{\delta_j^k|\xi|^2 -n\, \xi_j\xi_k}{|\xi|^{n+2}}  \hat{g}_{\varepsilon, j}(\xi)\overline{\hat{g}_{\varepsilon, k}(\xi)}\, d\xi\Bigr ),
\end{eqnarray}
where all integrals are absolutely convergent. 
By (\ref{NN}), 
$$\hat{\bf g}_\varepsilon(\xi) = \hat{\bf g}(\xi) -   e^{-|\xi|^2/2\varepsilon^2}\hat{\bf g}(0).$$
We note that for any $t>0$
$$
\int_{|\xi|> t}\ \frac{\delta_j^k|\xi|^2 -n\, \xi_j\xi_k}{|\xi|^{n+2}}   e^{-|\xi|^2/2\varepsilon^2}d\xi =0$$
and
$$
\int_{|\xi|> t}\ \frac{\delta_j^k|\xi|^2 -n\, \xi_j\xi_k}{|\xi|^{n+2}} e^{-|\xi|^2/\varepsilon^2}\, \hat{g}_j(0)\,  \overline{\hat{g}_k(0)}\,  d\xi =0.$$
Therefore
\begin{eqnarray*}
&&\int_{|\xi|> t}\ \frac{\delta_j^k|\xi|^2 -n\, \xi_j\xi_k}{|\xi|^{n+2}}\, \hat{g}_{\varepsilon,j}(\xi)\,  \overline{\hat{g}_{\varepsilon, k}(\xi)} \, d\xi\\
&& =\int_{|\xi|> t} \frac{\delta_j^k|\xi|^2 -n\, \xi_j\xi_k}{|\xi|^{n+2}}\, \hat{g}_{j}(\xi)\,  \overline{\hat{g}_{ k}(\xi)}\,  d\xi + O(\varepsilon)
\end{eqnarray*}

\noindent
uniformly with respect to $t$. Hence the value (\ref{2y}) tends to 
\begin{equation}\label{43w}
(2\pi)^{-n} \Bigl ( \sum _{1\leq j, k\leq n}\int \frac{\delta_j^k|\xi|^2 -n\, \xi_j\xi_k}{|\xi|^{n+2}}\, \hat{g}_{j}(\xi)\,  \overline{\hat{g}_{ k}(\xi)} d\xi \Bigr )
\end{equation}
as $\varepsilon \to 0_+$, where the integral is understood as the Cauchy  value.

\smallskip

Note that for $n>2$
$$
\bigl(|\xi|^2 - n\, \xi^2_k\bigr)|\xi|^{-2-n} = (2-n)^{-1}\frac{\partial ^2}{\partial\xi^2_k} \frac{1}{|\xi|^{n-2}} - \frac{|S^{n-1}|}{n}\delta(\xi)$$
and
$$\xi_j\xi_k\, |\xi|^{-2-n} = n^{-1} (n-2)^{-1} \frac{\partial ^2}{\partial\xi_j\, \partial \xi_k}\frac{1}{|\xi|^{n-2}}, \quad\quad {\rm for } \,\, j\neq k.
$$
Analogously, for $n=2$, 
\begin{equation}\label{4a}
\bigl(|\xi|^2 - 2\, \xi^2_k\bigr)|\xi|^{-4} = -\frac{\partial}{\partial\xi_k}\, \frac{\xi_k}{|\xi|^2} - \pi\, \delta(\xi),
\end{equation}
and
\begin{equation}\label{5a}
\xi_j\xi_k\, |\xi|^{-4} = -\frac{1}{2}\frac{\partial}{\partial\xi_j}\, \frac{\xi_k}{|\xi|^2} \qquad {\rm for } \,\, j\neq k.
\end{equation}

\noindent
Therefore,  in the case $n>2$,  we express (\ref{43w}) as
 
$$(2\pi)^{-n} \sum _{1\leq k\leq n}\int \left( \frac{1}{2-n}\, \frac{\partial ^2}{\partial \xi^2_k}\, \frac{1}{|\xi|^{n-2}} - \frac{|S^{n-1}|}{n}\, \delta(\xi)\right) {\hat{g_k}}(\xi)\, {\bar{\hat{g_k}}}(\xi)\, d\xi$$

$$-(2\pi)^{-n} n\sum_{j\neq k}\int \frac{1}{n(n-2)}\Bigl( \frac{\partial ^2}{\partial\xi_j\partial \xi_k}\, \frac{1}{|\xi|^{n-2}}\Bigr) {\hat {g_k}}(\xi){\bar{\hat{g_j}}}(\xi)\, d\xi.$$

\noindent
Using Parseval's formula once more, we write the limit of the right-hand side in  (\ref{2y}) as $\varepsilon\to 0_+$ in the form 

$$ \sum _{1\leq k\leq n}\int {g_k}(x)\, {\mathcal F}^{-1}_{\xi\to x}\left(\left( \frac{1}{2-n}\, \frac{\partial ^2}{\partial \xi^2_k}\, \frac{1}{|\xi|^{n-2}} - \frac{|S^{n-1}|}{n}\, \delta(\xi)\right) {\hat{g_k}}(\xi)\right) \, dx$$

\begin{equation}\label{6a}
-\sum_{j\neq k}\int \frac{1}{n-2}\,{g_j}(x)\, {\mathcal F}^{-1}_{\xi\to x}\left(\left( \frac{\partial ^2}{\partial\xi_j\partial \xi_k}\, \frac{1}{|\xi|^{n-2}}\right) {\hat {g_k}}(\xi)\right) \, dx,
\end{equation}
where ${\mathcal F}^{-1}$ means the inverse Fourier transform (see formula $(7.1.4)$ in \cite{[H]}). 
Since ${\mathcal F}^{-1} ({\hat u}\, {\hat v}) = u\ast v$, where $\ast$ denotes the convolution, we have
$$
{\mathcal F}^{-1}_{\xi\to x}\left(\left( \frac{\partial ^2}{\partial\xi_j\partial \xi_k}\, \frac{1}{|\xi|^{n-2}}\right) {\hat {h}}(\xi)\right) = -\left( x_j\, x_k\Bigl({\mathcal F}^{-1}_{\xi\to x} \frac{1}{|\xi|^{n-2}}\Bigr)\right)\ast h$$
\begin{equation}\label{7a}
 = -(2\pi)^{-n}(n-2) |S^{n-1}|\frac{x_j\, x_k}{|x|^2}\ast h
\end{equation}
for $1\leq j, k\leq n$.

\smallskip

Now let $n=2$. By (\ref{4a}) and Parseval's formula, we   present   (\ref{43w}) in the form analogous to (\ref{6a})
\begin{align}
& (2\pi)^{-2} \sum_{1\leq k\leq 2} \int \Bigl(-\frac{\partial}{\partial\xi_k}\, \frac{\xi_k}{|\xi|^2} - \pi\, \delta(\xi)\Bigr) {\hat{g_k}}(\xi)\, {\bar{\hat{g_k}}}(\xi)\, d\xi\nonumber\\
 &+   (2\pi)^{-2} \sum_{j\neq k} \int \Bigl(\frac{\partial}{\partial\xi_j}\, \frac{\xi_k}{|\xi|^2}\Bigr){\hat{g_k}}(\xi)\, {\bar{\hat{g_j}}}(\xi)\, d\xi\nonumber\\
&= \sum_{1\leq k\leq 2} \int {\mathcal F}_{\xi\to x}^{-1}\left(\bigl(-\frac{\partial}{\partial\xi_k}\, \frac{\xi_k}{|\xi|^2} - \pi\, \delta(\xi)\bigr) {\hat{g_k}}(\xi)\right) g_k(x)\, dx\nonumber\\
&+ \sum_{j\neq k} \int {\mathcal F}_{\xi\to x}^{-1}\left( \bigl(\frac{\partial}{\partial\xi_j}\, \frac{\xi_k}{|\xi|^2}\bigr){\hat{g_k}}(\xi)\right) g_j(x)\, dx.\label{26}
\end{align}

\noindent
We check directly that 

$${\mathcal F}_{\xi\to x}^{-1}\left( \bigl(\frac{\partial}{\partial\xi_j}\, \frac{\xi_k}{|\xi|^2}\bigr){\hat{h}}(\xi)\right) = 
i\, x_j {\mathcal F}_{\xi\to x}^{-1}\frac{\xi_k}{|\xi|^2}{\hat{h}}(\xi) = -(2\pi)^{-2}\frac{x_j\, x_k}{|x|^2}\ast h.$$

\noindent
Combining this with (\ref{7a}), we deduce from (\ref{6a}) and (\ref{26}) that for every $n\geq 2$ the limit of the expression (\ref{2y}) as $\varepsilon\to 0_+$ is equal to

\begin{eqnarray}\label{N0y}
&&\Bigl (\sum_{1\leq k \leq n}\int \left(\Bigl(\frac{|S^{n-1}|}{(2\pi)^n}\, \frac{x_k^2}{|x|^2} - \frac{|S^{n-1}|}{(2\pi)^nn}\Bigr)\ast g_k\right) g_k dx + \sum_{j\neq k} \frac{|S^{n-1}|}{(2\pi)^n}\Bigl(\frac{x_j\, x_k}{|x|^2}\ast g_k\Bigr)g_j dx \, \Bigr )\nonumber\\
\nonumber\\
&& = \frac{|S^{n-1}|}{(2\pi)^n}\Bigl (\int  \sum_{1\leq j, k\leq n} \Bigl(\frac{x_j\, x_k}{|x|^2}\ast g_j\Bigr)g_k\, dx - \frac{1}{n} \sum_{k=1}^n\Bigl(\int  g_k\, dx\Bigr)^2\Bigr) \nonumber\\
\nonumber\\
&&= \frac{|S^{n-1}|}{(2\pi)^n}\int \!\int \Bigl(  {\mathcal M}\Bigl(\frac{x-y}{|x-y|}\Bigr) {\bf g}(x), {\bf g}(y)\Bigr) dx\, dy,
\end{eqnarray}
where ${\mathcal M}(\omega)$ is the $(n\times n)$-matrix given by
\begin{equation}\label{n11}
{\mathcal M}(\omega) = (\omega_j\, \omega_k - n^{-1} \delta_j^k)_{j,k=1}^n.
\end{equation}
 Since the norm of ${\mathcal M}(\omega)$ does not exceed $(n-1)\, n^{-1}$, it follows that the absolute value of the last double integral is not greater than
\begin{equation}\label{9a}
\frac{|S^{n-1}|(n-1)}{(2\pi)^n\, n}\Bigl(\int |{\bf g}(x)|\, dx\Bigr)^2.
\end{equation} 
Hence (\ref{9a}) is a majorant for the left-hand side of (\ref{39w}). It remains to recall that $|S^{n-1}| = 2 \pi^{n/2}/\Gamma(n/2)$. 

\smallskip

$(ii)$ By (\ref{N1y}) and $(i)$,
\begin{eqnarray}\label{N3y}
\frac{|n-c_1|}{n} \mathop{\hbox {lim sup}}_{\varepsilon\to 0_+}\|{\bf g}_\varepsilon\|_{{\cal H}^{-n/2}}^2 &\leq& \frac{c_1}{n} \mathop{\hbox {lim sup}}_{\varepsilon\to 0_+} \Bigl | \|{\bf g}_\varepsilon\|_{{\cal H}^{-n/2}}^2 -n \|{\rm div}\, {\bf g}_\varepsilon\|_{{\cal H}^{-1-n/2}}^2\Bigr|\nonumber\\
&+& c_2\,  \|{\bf g}\|_{L^1}^2 \leq c_3\,  \|{\bf g}\|_{L^1}^2.
\end{eqnarray}
Since $L^1$ is not embedded into ${\cal H}^{-n/2}$, we have $c_1 =n$. 

\smallskip

 Suppose that (\ref{N1y}) holds. Then $c_1 =n$ and by (\ref{2y}) and (\ref{43w})
the inequality
\begin{equation}\label{N1w}
(2\pi)^{-n}\Bigl |\sum_{1\leq j,k\leq n}\int\ \frac{\delta_j^k|\xi|^2 -n\, \xi_j\xi_k}{|\xi|^{n+2}}  \hat{g}_j(\xi)\overline{\hat{g}_k(\xi)}\, d\xi\Bigr | \leq c_2\, \|{\bf g}\|_{L^1}^2
\end{equation}
holds for ${\bf g}\in C_0^\infty$ with the integral   understood as the Cauchy  value. It was shown in the proof of part $(i)$  that (\ref{43w}) is equal to (\ref{N0y}). Thus (\ref{N1w}) can be written as  the inequality
\begin{equation}\label{N4y}
\frac{|S^{n-1}|}{(2\pi)^n}\Bigl |\int \!\int \Bigl(  {\mathcal M}\Bigl(\frac{x-y}{|x-y|}\Bigr) {\bf g}(x), {\bf g}(y)\Bigr) dx\, dy\Bigr | \leq c_2\,  \|{\bf g}\|_{L^1}^2, 
\end{equation}

\noindent
where the matrix ${\mathcal M}$ is defined by (\ref{n11}). 
Let $\theta$ denote the north pole of $S^{n-1}$, i.e. $\theta = (0, \ldots, 0, 1)$. We choose the vector function ${\bf g}$ in (\ref{N4y}) as $(0, \ldots, \eta(|x|)\varphi(x/|x|))$, where $\eta\in C_0^\infty([0,\infty))$, $\eta\geq 0$, and $\varphi$ is a regularization of the $\delta$-function on $S^{n-1}$ concentrated at $\theta$. Then (\ref{N4y}) implies
\begin{eqnarray*}
&&\frac{|S^{n-1}|}{(2\pi)^n}\Bigl |\int_0^\infty\! \int_0^\infty m_{nn}\Bigl(\frac{\rho -r}{|\rho -r|}\theta\Bigr) \eta(r) \; r^{n-1}\, \eta(\rho)\, \rho^{n-1} drd\rho \Bigr |\\
\\
&&\leq c_2 \Bigl | \int_0^\infty \eta (r) r^{n-1} dt \Bigr |^2, 
\end{eqnarray*}
and since $m_{nn} (\pm \theta) = 1-1/n$, we obtain  $c_2\geq (1- 1/n)(2\pi)^{-n} |S^{n-1}|$.  

\smallskip

$(iii)$ By (\ref{Nd}), we change $n^{-1}$ in (\ref{N0y}) for $1/2$ and notice that the norm of the matrix $(\omega_j\omega_k -\delta_j^k/2)_{j,k =1}^n$ equals $1/2$. Inequality (\ref{46a}) follows. $\square$

\medskip

As an immediate consequence of Theorem 3 $(iii)$, we derive

\smallskip

{\bf Corollary.} {\it Let $u$ be a scalar function in $C_0^\infty$. Then}
\begin{equation}\label{gu}
\|u\|_{{\cal H}^{1-n/2}} \leq  \Bigl( \frac{(2\sqrt{\pi})^{-n}}{\Gamma(n/2) (n-1)}\Bigr)^{1/2} \|\nabla u\|_{L^1}.
\end{equation}

\smallskip

{\bf Proof.} It suffices to put ${\bf g} = \nabla u$ in (\ref{46a}) and note that
$$\|{\bf g} \|_{{\cal H}^{-n/2}} = \|\nabla u\|_{{\cal H}^{-n/2}} = \|u\|_{{\cal H}^{1-n/2}}$$
and

$$\|{\rm div}\, {\bf g}\|_{{\cal H}^{-1-n/2}} = \|\Delta u\|_{{\cal H}^{-1-n/2}} = \|u\|_{{\cal H}^{1-n/2}}.\qquad\qquad \quad\square $$

\medskip

{\bf Remark 2.}  Passing from quadratic to sesquilinear forms in  the proof of Theorem 3  $(i)$ leads to the identity
\begin{equation}\label{cr}
(-\Delta)^{-n/2} \bigl ({\bf g} + n(-\Delta)^{-1} \nabla {\rm div}\, {\bf g}\bigr ) (x) = \frac{2^{1-n}\pi^{-n/2}}{\Gamma(n/2)} \int {\mathcal N}\Bigl(\frac{x-y}{|x-y|}\Bigr) {\bf g}(y)\, dy 
\end{equation}
for all ${\bf g}\in C_0^\infty$ orthogonal to $1$. The kernel ${\mathcal N} (\omega)$ is the matrix function $(\omega_j\omega_k)_{j,k =1}^n$. Needless to say, if additionally ${\bf g}$ is divergence free, we have the representation
\begin{equation}\label{cru}
(-\Delta)^{-n/2} {\bf g} (x) = \frac{2^{1-n}\pi^{-n/2}}{\Gamma(n/2)} \int {\mathcal N}\Bigl(\frac{x-y}{|x-y|}\Bigr) {\bf g}(y)\, dy. 
\end{equation}

\smallskip

Another consequence of the  identity (\ref{cr}) is obtained by putting   ${\bf g} = \nabla u$ in it, where $u$ is a scalar function in $C_0^\infty$. Then 
\begin{equation}\label{cre}
(-\Delta)^{-n/2} \nabla u(x) =  \frac{2^{1-n}\pi^{-n/2}}{(1-n) \Gamma(n/2)} \int {\mathcal N}\Bigl(\frac{x-y}{|x-y|}\Bigr) \nabla u(y)\, dy. 
\end{equation}

\medskip

{\bf Remark 3.}  If ${\rm div}\, {\bf g}\in {\cal H}^{-1-n/2}$ and ${\bf g}\in L^1$, then ${\bf g}$ is orthogonal to one (see the beginning of the proof of Theorem 2). On the other hand, 
 even if ${\bf g}\in C_0^\infty$ but 
$$\int {\bf g}(y)\, dy \neq 0,$$
 both norms $\|{\rm div}\, {\bf g}\|_{{\cal H}^{-1-n/2}}$ and $\|{\bf g}\|_{{\cal H}^{-n/2}}$ are infinite. The estimate (\ref{39w}) shows that the formal expression

\begin{equation}\label{Nye}
\|{\bf g}\|_{{\cal H}^{-n/2}}^2 -n \|{\rm div}\, {\bf g}\|_{{\cal H}^{-1-n/2}}^2
\end{equation}

\noindent
can be given sense as the finite limit $\varepsilon\to 0_+$on the left-hand side of (\ref{39w}). One can see that the limit   does not change if (\ref{NN}) is replaced by 
$$ {\bf g}_\varepsilon(x) = {\bf g}(x) -  \varepsilon^n \, \eta(\varepsilon x) \int {\bf g}(y) \, dy,$$
where $\eta$ is an arbitrary function in the Schwartz space ${\cal S}$ normalized by 
$$\int \eta(y)\, dy =1.$$

\smallskip

By Theorem 3 $(iii)$ and a 
 duality argument, similar to that used in \cite{[BB3]},  one can arrive at the following existence result which is supplied with a  proof for reader's convenience.  

\medskip

{\bf Proposition 3.}  {\it For any vector function ${\bf u} \in {\cal H} ^{n/2}$ there exists a vector function  ${\bf v}\in L^\infty$ and a scalar function $\varphi\in {\cal H}^{1+ n/2}$ satisfying} ${\bf u} = {\bf v} + {\rm grad}\, \varphi$. 

\smallskip
{\bf Proof.} By ${\cal B}$ we denote the Banach space of the pairs $\{{\bf g}, k\}\in L^1\times {\cal H}^{-1-n/2}$ endowed with the norm
$$\| \{{\bf g}, k\}\|_{{\cal B}} = \|{\bf g}\|_{L^1} + \|k\|_{{\cal H}^{-1-n/2}}.$$
Representing $\{{\bf g}, k\}$ as $\{{\bf g}, 0\} + \{0, k\}$, we see that an arbitrary linear functional on ${\cal B}$ can be given by
\begin{equation}\label{oi}
\int ({\bf v}, \, {\bf g}) \, dx + \int \varphi\, k\, dx,
\end{equation}
where ${\bf v}\in L^\infty$ and $\varphi\in {\cal H}^{1+n/2}$.  The range of the operator
$$L^1\cap {\cal H}^{-n/2} \ni {\bf g} \to \{{\bf g}, - {\rm div}\, {\bf g} \},$$
which is a closed subspace  of ${\cal B}$,  will be denoted by $S$.

\smallskip

Any vector-valued  function ${\bf u}\in {\cal H}^{n/2}$ generates the continuous functional
\begin{equation}\label{oi1}
f({\bf g}) = \int ({\bf u}, \, {\bf g}) \, dx 
\end{equation}
on the space ${\cal H}^{-n/2}$.  By (\ref{46a}),
\begin{equation}\label{oi2}
| f({\bf g}) | \leq c_n \, \|{\bf u} \|_{{\cal H}^{n/2}}\bigl (\|{\bf g}\|_{L^1} + \|{\rm div}\, {\bf g}\|_{{\cal H}^{-1-n/2}}\bigr).
\end{equation}
We introduce the functional $\Phi$ by 
$$\Phi(\{{\bf g}, k\}) := f({\bf g}) \quad {\rm for} \quad k= -{\rm div}\, {\bf g},$$
i.e. $\Phi$ is defined on   $S$.  Being prescribed on a closed subspace of ${\cal B}$, this functional is bounded in the norm of ${\cal B}$ because of (\ref{oi2}). By the Hahn-Banach theorem, $\Phi$ can be  extended with preservation of the norm  onto the whole space ${\cal B}$.  Using (\ref{oi}) and (\ref{oi1}), we see that there exist ${\bf v}\in L^\infty$ and $\varphi\in {\cal H}^{1+n/2}$ such that,   for
 all ${\bf g}\in L^1 \cap {\cal H}^{-n/2}$,
 $$\int ({\bf u}, \, {\bf g}) \, dx = \int \bigl ( ({\bf v}, \, {\bf g})  - \varphi\, {\rm div}\, {\bf g} \bigr) dx.$$
 The result follows. $\square$

\medskip

The next  assertion    guarantees the existence of the solution ${\bf u}\in {\cal H}^{2-n/2}$ to the equation $ -\Delta {\bf u} = {\bf f}$ provided that ${\bf f}$ 
 {\it is a vector field in $L^1$ subject to} 
  ${\rm div}\, {\bf f}\in {\cal H}^{-1-n/2}$.

\smallskip

{\bf Proposition 4.}  {\it Under the condition on ${\bf f}$ just mentioned,  the inequality
\begin{equation}\label{Nt}
\Bigl | \| (-\Delta)^{-1} {\bf f} \|^2_{{\cal H}^{2-n/2}} -n \|{\rm div}\, {\bf f}\|^2_{{\cal H}^{-1-n/2}} \Bigr | \leq \frac{ (2\sqrt{\pi})^{-n}}{\Gamma(n/2)}\|{\bf f} \|_{L^1}^2
\end{equation}
holds.}

\smallskip

{\bf Proof.} It suffices to replace ${\bf g}$ by ${\bf f}$ in (\ref{46a}).  
$\square$

\smallskip

In the forthcoming  Theorem $4$ we obtain an estimate   which leads by duality to the following existence result. Its proof is quite similar to that of Proposition  3  and is omitted.

\medskip

{\bf Proposition 5.}  {\it Let ${\bf f}$ be a divergence free vector function in $\Bbb{R}^3$ from the space ${\cal H}^{1/2}$. Then the equation 
$${\rm curl}\, {\bf u} = {\bf f}\quad {\rm in} \,\,\, \Bbb{R}^3$$
has a solution in } ${\cal H}^{3/2}\cap L^\infty$. 

\medskip

{\bf Theorem 4.} {\it Let
\begin{equation}\label{Vy}
{\rm curl}\, {\bf w} = {\bf f} + {\bf g} \quad {\rm in} \,\,\, \Bbb{R}^3,
\end{equation}
where 
$${\rm div}\, {\bf w} =0, \quad  {\bf f}\in {\cal H}^{-3/2}(\Bbb{R}^3)$$
 and
 $$ {\bf g}\in L^1(\Bbb{R}^3), \quad \int_{\Bbb{R}^3} {\bf g}(y)dy =0.$$
 Then}
 
\begin{equation}\label{Wy}
\Bigl |\|\Delta{\bf w} + {\rm curl}\, {\bf f} \|_{{\cal H}^{-5/2}}^2 - 2\, \|{\rm div}\, {\bf f}\|_{{\cal H}^{-5/2}}^2\Bigr| \leq \frac{1}{4\pi^2}\Bigl(\int_{\Bbb{R}^3}|{\bf g}(x)|dx\Bigr)^2.
\end{equation}

\smallskip

{\bf Proof.} Since ${\rm curl}^2{\bf w} = -\Delta {\bf w}$, we have by (\ref{Vy}) that $-\Delta {\bf w} = {\rm curl}\, {\bf f} + {\rm curl}\, {\bf g}$. Using the identity ${\rm div}\, {\rm curl}\, {\bf w} =0$, we see that ${\rm div}\, {\bf f} + {\rm div}\, {\bf g} =0$. Therefore,
$$
\|\Delta{\bf w} + {\rm curl}\, {\bf f} \|_{{\cal H}^{-5/2}}^2 - 2\, \|{\rm div}\, {\bf f}\|_{{\cal H}^{-5/2}}^2 =\| {\rm curl}\, {\bf g} \|_{{\cal H}^{-5/2}}^2 - 2\, \|{\rm div}\, {\bf g}\|_{{\cal H}^{-5/2}}^2.$$

\noindent
The right-hand side can be written in the form
$$(2\pi)^{-3} \Bigl | \int_{\Bbb{R}^3}\Bigl (|\xi\times \hat{\bf g} |^2 - 2\, |(\xi, \hat{\bf g})|^2\Bigr)\frac{d\xi}{|\xi|^5} \Bigr |= (2\pi)^{-3} \Bigl | \int_{\Bbb{R}^3}  \Bigl (|\xi|^2 |\hat{\bf g} |^2 - 3\, |(\xi, \hat{\bf g})|^2\Bigr)\frac{d\xi}{|\xi|^5}\Bigr |.$$ 
This value is a particular case of (\ref{2y}) for $n=3$ and hence it does not exceed
$$\frac{1}{4\pi^2}\Bigl(\int_{\Bbb{R}^3}  |{\bf g}(x)| \, dx\Bigr)^2$$
(see  the proof of Theorem 3 $(iii)$). $\square$

\medskip

{\bf Remark 5.} It is natural to ask how the results of the present section change if the role of the homogeneous space ${\cal H}^l$ is played by the standard Sobolev space $H^l$ endowed with the norm
$$\|\phi\|_{H^l} : = \Bigl( \int |\hat{\phi}(\xi)|^2 (|\xi|^2 +1)^{l/2} d\xi \Bigr)^{1/2}.$$

Restricting ourselves to Theorem 3,  we check directly that
\begin{eqnarray*}
&&\Bigl | \lim\limits_{\varepsilon\to 0_+}\bigl(\|{\bf g}_\varepsilon\|^2_{{ H}^{-n/2}} - n\, \|{\rm div} \, {\bf g}_\varepsilon \|^2_{{ H}^{-1-n/2}} \bigr) \Bigr | \\
\\
&&= (2\pi)^{-n}\Bigl |\sum_{1\leq j,k\leq n}\int \frac{\delta_j^k(|\xi|^2+1) -n\, \xi_j\xi_k}{(|\xi|^{2}+1)^{1+n/2}}  \hat{g}_{ j}(\xi)\overline{\hat{g}_{ k}(\xi)}\, d\xi\Bigr |,
\end{eqnarray*}
which in its turn is equal to
\begin{eqnarray*}
&&(2\pi)^{-n} (n-2)^{-1} \Bigl |\sum_{1\leq j,k\leq n}\int \frac{\partial ^2}{\partial\xi_j \partial\xi _k}(|\xi|^{2}+1)^{(2-n)/2}\hat{g}_{ j}(\xi)\overline{\hat{g}_{ k}(\xi)}\, d\xi\Bigr |\\
\\
&&= c\Bigl |\int \sum_{1\leq j,k\leq n} \frac{x_j - y_j}{|x-y|}\frac{x_k - y_k}{|x-y|} |x-y| K_1(|x-y|) g_j(x) g_k(y) \, dxdy\Bigr |,
\end{eqnarray*}
where $K_1$ is the modified Bessel function of the third kind.  Since the function $t\, K_1(t)$ is bounded, we obtain
$$\Bigl | \lim\limits_{\varepsilon\to 0_+}\bigl(\|{\bf g}_\varepsilon\|^2_{{ H}^{-n/2}} - n\, \|{\rm div} \, {\bf g}_\varepsilon \|^2_{{ H}^{-1-n/2}} \bigr) \Bigr |  \leq c(n) \Bigl( \int |{\bf g}(x)|\, dx \Bigr)^2.$$
Needless to say, this inequality becomes
$$\Bigr |\|{\bf g}\|^2_{{ H}^{-n/2}} - n\, \|{\rm div} \, {\bf g} \|^2_{{ H}^{-1-n/2}}  \Bigr |  \leq c(n) \Bigl( \int |{\bf g}(x)|\, dx \Bigr)^2$$
if 
 the last norm of ${\rm div}\, {\bf g}$ is finite.

\bigskip

\section*{Acknowledgement}

I cordially thank  Ha\"im Brezis whose questions infused me with the idea to write this article. I gratefully acknowledge referee's poignant comments.

\end{document}